\documentclass[11pt,a4paper]{article}
\usepackage{amssymb}
\title{Chow's Theorem for Linear Spaces}
\author{Hans Havlicek}
\date{}

\newtheorem{theo}{Theorem}
\newtheorem{lemma}{Lemma}

\newcommand{\Ecal}{{\cal E}}

\newcommand{\Lcal}{{\cal L}}
\newcommand{\Mcal}{{\cal M}}
\newcommand{\Ncal}{{\cal N}}

\newcommand{\Pcal}{{\cal P}}

\newcommand{\Scal}{{\cal S}}

\newcommand{\Xcal}{{\cal X}}

\newcommand{\rel}{\sim}
\newcommand{\rrel}{\approx}

\newcommand{\PL}{\mbox{$(\Pcal,\Lcal)$}}
\newcommand{\PLL}{\mbox{$(\Pcal',\Lcal')$}}%

\newcommand\inv{^{-1}}

\newcommand{\abb}[3]{\mbox{$#1\,:\,#2\rightarrow#3$}}
\newcommand{\Abb}[5]{\mbox{$#1\,:\,#2\rightarrow#3,\;#4\mapsto #5$}}
\newcommand{\proof}{{\em Proof. }}
\newcommand{\proofend}{$\Box$}
\newcommand{\zitat}[4]{\bibitem{#1}{\sc #2}: {\sl #3\/}.
#4.\vspace{-0.7em}}

\begin{document}
\maketitle

   \begin{abstract}\sloppy
   If \abb{\varphi}{\Lcal}{\Lcal'} is a bijection from the set of
   lines of a linear space \PL\ onto the set of lines of a linear
   space \PLL\ ($\dim\PL$,  $\dim\PLL\geq 3$), such that intersecting
   lines go over to intersecting lines in both directions, then
   $\varphi$ is arising from a collineation of \PL\ onto \PLL\ or a
   collineation of \PL\ onto the dual linear space of \PLL. However,
   the second possibility can only occur when \PL\ and \PLL\ are
   $3$--dimensional generalized projective spaces.
   \end{abstract}


\section{Introduction}\label{Intro}

The Grassmannian $\Gamma_{k,n}$ formed by the $k$--dimensional
subspaces of an $n$--dimen\-sional projective space ($1\leq k\leq
n-2)$ carries the structure of a metric space: Two $k$--subspaces are
at distance $d$ if their meet is $(k-d)$--dimensional \cite[p.\
807]{He95}; subspaces with distance $1$ are called {\em adjacent}. Any
collineation and, however only if $n=2k+1$, any duality of the
underlying projective space yields a bijection on $\Gamma_{k,n}$ which
is distance preserving in both directions. Conversely, any bijection
of $\Gamma_{k,n}$ which is adjacency preserving in both directions
arises in this way. The last result is due to {\sc W.--L.\ Chow}; cf.
\cite{Ch49} or \cite[80--82]{Di71}. Thus Chow's theorem may be seen as
an early result in a discipline which now is called {\em
characterizations of geometrical transformations under mild
hypotheses} \cite{Be92}, \cite{Be94}, \cite{Le95}. Alternatively, one
may consider $\Gamma_{k,n}$ as the set of vertices of a graph where
two vertices are joined by an edge if, and only if, they represent
adjacent subspaces. However, we shall not adopt this point of view.

Other results in the spirit of Chow's theorem can be found in
\cite{Be92}, \cite{Br88}, \cite{Ch49}, \cite[82--88]{Di71},
\cite{Ha95}, and \cite{Ha96}.

In an arbitrary linear space there are several concepts of
``dimension''. Cf.\ the remarks in \cite[p.\ 73]{Bu95}. Hence one may
ask if Chow's theorem still holds true for linear spaces. We shall not
deal with this question in full generality, as we focus our attention
on the case of $1$--dimensional subspaces, i.e., the lines of a linear
space. Loosely speaking, here the answer is in the affirmative. We
refer to Theorem \ref{KAPPA} for a precise statement.

\section{Preliminaries}\label{Basic}

Let \PL\ be a linear space. We shall stick to the terminology but not
to the notation used in \cite{KSW73}: Given a subset
$\Scal\subset\Pcal$ then $\langle\Scal\rangle$ denotes the {\em span}
of $\Scal$, i.e., the intersection of all subspaces containing
$\Scal$. The {\em join} of subsets $\Scal_1,\Scal_2\subset\Pcal$ is
defined as $\langle\Scal_1\cup\Scal_2\rangle=:\Scal_1\vee\Scal_2$. The
{\em dimension} of \PL\ is given by%
   \footnote{Cf., however, the definition in \cite[p.\ 9]{BB93}.}
   \begin{equation}
   \dim \PL := \min\{\#\Xcal - 1\mid \Xcal\subset\Pcal,\;
   \langle\Xcal\rangle=\Pcal\}.
   \end{equation}
A {\em plane} of \PL\ is a two--dimensional subspace.

A linear space \PL\ is called an {\em exchange space}, if it satisfies
the exchange axiom: Given any subset $\Scal\subset\Pcal$ and points
$A,B\in\Pcal$ then
   \begin{equation}\label{AUSTAUSCHLEMMA}
   B\in(\Scal\vee \{A\})\setminus \langle\Scal\rangle \Longrightarrow
   A\in\Scal\vee \{B\}.
   \end{equation}
We observe that (\ref{AUSTAUSCHLEMMA}) is true in any linear space, if
$A,B$ are collinear with some point $X\in\langle\Scal\rangle$.

A 2--dimensional linear space such that any two distinct lines have at
least one point in common is called a {\em generalized projective
plane}. A {\em generalized projective space} is a linear space \PL\
where any plane is a generalized projective plane \cite[p.\ 9]{BB93};
cf.\ also \cite[p.\ 39]{BC95}. The following property of generalized
projective spaces is immediate: Given a subset $\Scal\neq\emptyset$
and a point $X\in\Pcal\setminus\langle\Scal\rangle$ then
   \begin{equation}\label{VERBIND}
   \Scal\vee\{X\}=\bigcup_{Y\in\langle\Scal\rangle}(\{X\}\vee \{Y\}).
   \end{equation}
Consequently, \PL\ is an exchange space.

We shall encounter 3--dimensional generalized projective spaces in
Theorem \ref{KAPPA}. If \PL\ is such a space,  then any plane
$\Ecal\subset\Pcal$ and any line $y\in\Lcal$ have a common point
according to (\ref{VERBIND}). Moreover, any two distinct planes meet
at a unique line. Conversely, given a line $a\in\Lcal$ then denote by
$a^\ast$ the {\em pencil of planes} with axis $a$, i.e., the set of
all planes passing through $a$. As \PL\ is an exchange space, any two
distinct points on $a$ extend to a basis of \PL, whence $\#a^\ast\geq
2$. Thus it is possible to define the {\em dual linear space} of \PL\
as follows: Let $\Pcal^\ast$ be the set of all planes of \PL\ and let
$\Lcal^\ast$ be the set of all pencils of planes \cite[p.\ 72]{Bu95}.

\section{Line Geometry}\label{LINE}

Let \PL\ be a linear space. Lines $a,b\in \Lcal$ are called {\em
related} ($a\rel b$) if $a\cap b\neq\emptyset$. Distinct related lines
$a,b$ are denoted by $a\rrel b$ and will be named {\em adjacent}. The
pair $(\Lcal,\rel)$ is a {\em Pl\"ucker space}%
   \footnote{We refrain from assuming that $\Lcal$ is non--empty.}
in the sense of {\sc W.\ Benz} \cite[p.\ 199]{Be92}.

Let $\Ncal\subset \Lcal$ be a set of mutually related lines; we shall
use the term {\em related set} for such an $\Ncal$. It is easily seen,
by virtue of Zorn's lemma, that there is a maximal related set
$\Mcal\subset\Lcal$ containing $\Ncal$. If $-1\leq\dim\PL\leq 1$, then
$\Lcal$ is the only maximal related set. With $A\in\Pcal$ write
$\Lcal(A)$ for the {\em star of lines} with vertex $A$, i.e., the set
of all lines of $\Lcal$ running through $A$. Any star of lines is a
maximal related set for $\dim\PL\neq 2$; in the 2--dimensional case a
star of lines obviously is a related set. It may be maximal (e.g.\ in
affine planes) or not maximal (e.g.\ in projective planes). If
$\dim\PL\geq 2$, then any trilateral is contained in a maximal related
set which cannot be a star of lines.
\begin{lemma}\label{MAXI}
   Let $\dim\PL\geq 2$ and let $\Mcal\subset\Lcal$ be a maximal
   related set different from a star of lines. Then the subspace
   $\Ecal$ spanned by the lines of $\Mcal$ is a plane satisfying
      \begin{equation}\label{MAXI1}
      \Ecal = x\vee y \mbox{ for all } x,y\in\Mcal,\; x\neq y.
      \end{equation}
   \end{lemma}
\proof  Choose any line $w\in\Mcal\setminus\{x,y\}$. If $\{x,y,w\}$ is
a trilateral, then $w\subset x\vee y\subset \Ecal$. Otherwise, there
exists a line $z\in\Mcal\setminus\{x,y\}$ such that $\{x,y,z\}$ is a
trilateral, since $\Mcal$ cannot be contained in a star of lines.
Hence $w\subset x\vee z = x\vee y\subset \Ecal$ and
   \begin{equation}
   \Ecal=\bigvee_{w\in\Mcal} w\subset x\vee y\subset\Ecal
   \end{equation}
is a plane.
\proofend

\vspace{1em}

For $\dim\PL\geq 3$ the maximal related sets in $\Lcal$ fall into two
classes: stars of lines and the {\em coplanar maximal related sets}
described in Lemma \ref{MAXI}. This result is crucial in proving
Theorem \ref{KAPPA}, since stars allow to recover the points of \PL\
in terms of line geometry. If $\dim\PL=2$, then this distinction of
the maximal related sets is in general not available. Hence
$2$--dimensional linear spaces will not appear in the next section.

\section{Chow's Theorem}\label{ISOM}

Let \PL\ and \PLL\ be linear spaces. Any collineation
$\abb\kappa{\Pcal}{\Pcal'}$ gives rise to a bijection
   \begin{equation}
   \Abb\varphi{\Lcal}{\Lcal'}{\{A\}\vee \{B\}}{\{A^\kappa\}\vee
   \{B^\kappa\}} \; (A,B\in\Pcal,\; A\neq B)
   \end{equation}
taking related lines to related lines in both directions. We shall
prove the following converse:
   \begin{theo}\label{KAPPA}
   Let \PL\ and \PLL\ be linear spaces such that  $\dim\PL\geq 3$ and
   $\dim\PLL\geq 3$. Suppose that $\abb\varphi\Lcal{\Lcal'}$ is a
   bijection with the property
      \begin{equation}
      a\rel b \Longleftrightarrow
      a^\varphi\rel' b^\varphi \mbox{ for all } a,b\in\Lcal.
      \end{equation}
   Then the following assertions hold true:
      \begin{enumerate}
      \item Under $\varphi,\varphi\inv$ maximal related sets go over
      to maximal related sets.
      \item If $\varphi$ maps one star of lines onto a star of lines,
      then%
      \footnote{Whenever it is convenient, we do not distinguish
      between a point $X$ and the set $\{X\}$.}
      \begin{equation}\label{DEFKAPPA}
         \Abb\kappa\Pcal{\Pcal'}{a\cap b}{a^\varphi\cap b^\varphi}\;
         (a,b\in\Lcal,\; a\rrel b)
      \end{equation}
      is a collineation.
      \item If $\varphi$ maps one star of lines onto a coplanar
      maximal related set, then \PLL\ is a $3$--dimensional
      generalized projective space. The mapping
      \begin{equation}\label{DEFDELTA}
         \Abb\delta\Pcal{\Pcal'^\ast}{a\cap b}{a^\varphi\vee
         b^\varphi}\; (a,b\in\Lcal,\; a\rrel b)
      \end{equation}
      is a collineation onto the dual linear space of \PLL. Therefore
      \PL\ is a $3$--dimensional generalized projective space too and
      $\delta$ is a correlation of \PL\ onto \PLL.
   \end{enumerate}
   \end{theo}
\proof {\em Ad 1.}
This is obviously true.

{\em Ad 2.}
Suppose that there is a point $A\in\Pcal$ with $(\Lcal(A))^\varphi$
being a star of lines, say $\Lcal'(A')$ where $A'\in\Pcal'$.

Choose any point $B\in\Pcal\setminus\{A\}$. Assume that
$(\Lcal(B))^\varphi$ is not a star of lines. By Lemma \ref{MAXI},
$(\Lcal(B))^\varphi$ is a coplanar maximal related set. Write $\Ecal'$
for the plane spanned by the lines of $(\Lcal(B))^\varphi$. We infer
from $(\Lcal(A)\cap\Lcal(B))^\varphi = (\{A\}\vee \{B\})^\varphi$ that
$A'\in\Ecal'$. Since $\dim\PLL\geq 3$, there exists a point
$X'\in\Pcal'\setminus\Ecal'$. Choose $x\in\Lcal(A)$ such that
$x^\varphi=\{A'\}\vee \{X'\}$. Then $x\neq \{A\}\vee \{B\}$. There is
a point $Y\in x\setminus\{A\}$. Put $y:=\{Y\}\vee
\{B\}\notin\Lcal(A)$. Therefore $y^\varphi\notin\Lcal'(A')$ and
$A'\notin y^\varphi$. Hence
   \begin{equation}
   x^\varphi\cap y^\varphi= x^\varphi\cap\Ecal'\cap y^\varphi=
   \{A'\}\cap y^\varphi = \emptyset.
   \end{equation}
This implies $x^\varphi\not\rel' y^\varphi$ in contradiction to $x\rel
y$.

Thus we have established that $(\Lcal(B))^\varphi$ is a star of lines.
Conversely, given a point $B'\in\Pcal'\setminus\{A'\}$ it follows in
the same manner that $(\Lcal'(B'))^{\varphi\inv}$ is a star of lines.
The previous discussion shows that (\ref{DEFKAPPA}) is a well--defined
surjection. If points $Q,R\in\Pcal$ are distinct, then
   \begin{equation}
   \#(\Lcal(Q)\cap\Lcal(R))=\#((\Lcal(Q))^\varphi\cap(\Lcal(R))^\varphi)=1,
   \end{equation}
whence $Q^\kappa\neq R^\kappa$. Three mutually distinct points
$Q,R,S\in\Pcal$ are collinear if, and only if,
$\#(\Lcal(Q)\cap\Lcal(R)\cap\Lcal(S))=1$. This in turn is equivalent
to the collinearity of $Q^\kappa,R^\kappa,S^\kappa\in\Pcal'$. Hence
$\kappa$ is a collineation.

{\em Ad 3.}
{\em (a)}
Suppose that there is a point in $A\in\Pcal$ with $(\Lcal(A))^\varphi$
being a coplanar maximal related set. By the second part of the
present proof, $\varphi$ maps all stars in \PL\ onto coplanar maximal
related sets in \PLL. We infer from (\ref{MAXI1}) that
(\ref{DEFDELTA}) is a well--defined mapping.

{\em (b)}
Let $\Ecal'$ be any plane of \PLL. There exist adjacent lines
$a',b'\in\Lcal'$ with $a'\vee b'=\Ecal'$. By putting
   \begin{equation}\label{MAB}
   \{E\}:= a'^{\varphi\inv}\cap b'^{\varphi\inv}\mbox{ and }
   \Mcal'(a',b'):= (\Lcal(E))^\varphi,
   \end{equation}
we obtain a point in $\Pcal$ and a coplanar maximal related set
spanning $\Ecal'$, respectively. We deduce from (\ref{DEFDELTA}) that
$E^\delta=\Ecal'$, whence $\delta$ is surjective.

{\em (c)}
Choose any line $y'\in\Lcal'$. Then there is a line $x\in\Lcal(E)$
with $x\rel y'^{\varphi\inv}$. Hence any $y'\in\Lcal'$ is related to
at least one line of $\Mcal'(a',b')$.

{\em (d)}
There exists a point $Y'\in\Pcal'\setminus\Ecal'$. By (c), each line
on $Y'$ meets $\Ecal'$ at some point, whence $\dim\PLL=3$.

{\em (e)}
By (c), each line meets every plane. We infer from (d) that no plane
is properly contained in $\Ecal$, whence each line $c'\subset\Ecal'$
can be written as $c'=\Ecal'\cap (c'\vee Y')$. So any two distinct
lines of $\Ecal'$ have a common point and hence $\PLL$ is a
$3$--dimensional generalized projective space.

{\em (f)}
Let $P,Q\in\Pcal$ be distinct. By (e), $(\Lcal(P))^\varphi$ equals the
set of all lines within the generalized projective plane $P^\delta$.
The bijectivity of $\varphi$ forces that $P^\delta\cap
Q^\delta=(\{P\}\vee \{Q\})^\varphi$ is a single line. Thus
$P^\delta\neq Q^\delta$. This means that $\delta$ is injective.

{\em (g)}
Three mutually distinct points $Q,R,S\in\Pcal$ are collinear if, and
only if, $\#(\Lcal(Q)\cap\Lcal(R)\cap\Lcal(S))=1$. This in turn is
equivalent to the ``collinearity'' of
$Q^\delta,R^\delta,S^\delta\in\Pcal'^\ast$ with respect to the dual
linear space of \PLL. Hence $\delta$ is a collineation of \PL\ onto
$(\Pcal'^\ast,\Lcal'^\ast)$. Finally, the dual linear space of \PLL\
is a $3$--dimensional generalized projective space too, whence the
assertion on \PL\ follows.
\proofend

\vspace {1em}

\noindent
Hans Havlicek,
Abteilung f\"ur Lineare Algebra und Geometrie,
Technische Universit\"at,
Wiedner Hauptstra{\ss}e 8--10/1133,
A--1040 Wien, Austria\\
EMAIL: {\tt havlicek@geometrie.tuwien.ac.at}

\end{document}